\documentclass[12pt]{article}
\usepackage{amsmath, amssymb, a4wide, color}

\author{Andreas H. Hamel \footnote{Yeshiva University New York, Department of Mathematical
Sciences, Belfer Hall, 2495 Amsterdam Avenue, New York, NY 10033, hamel@yu.edu, + 1 212 960 5400 ext. 6919,}
 \and Andreas L\"ohne  \footnote{Martin-Luther-University Halle-Wittenberg, Department of Mathematics,
andreas.loehne@mathematik.uni-halle.de}}

\title{Lagrange Duality in Set Optimization}

\date{}

\newtheorem{theorem}{Theorem}

\newtheorem{remark}[theorem]{Remark}
\newtheorem{lemma}[theorem]{Lemma}
\newtheorem{definition}[theorem]{Definition}
\newtheorem{proposition}[theorem]{Proposition}
\newtheorem{example}[theorem]{Example}

\numberwithin{equation}{section} \numberwithin{theorem}{section}

\newcommand{\of}[1]{\ensuremath{\left( #1 \right)}}

\newcommand{\cb}[1]{\ensuremath{ \left\{ #1 \right\} }}
\newcommand{\sqb}[1]{\ensuremath{ \left[ #1 \right] }}

\newcommand{\bs}{\backslash}

\newcommand{\pend}{\hfill $\square$}

\newcommand{\vp}{\ensuremath{\varphi}}

\newcommand{\R}{\mathrm{I\negthinspace R}}
\newcommand{\OLR}{\overline{\mathrm{I\negthinspace R}}}
\newcommand{\N}{\mathrm{I\negthinspace N}}

\newcommand{\dom}{{\rm dom \,}}

\newcommand{\gr}{{\rm graph \,}}
\newcommand{\cl}{{\rm cl \,}}
\newcommand{\Int}{{\rm int\,}}

\newcommand{\co}{{\rm co \,}}

\newcommand{\Min}{{\rm Min \,}} 
\newcommand{\Max}{{\rm Max \,}}

\newcommand{\F}{\ensuremath{\mathcal{F}}}
\newcommand{\G}{\ensuremath{\mathcal{G}}}
\renewcommand{\P}{\ensuremath{\mathcal{P}}}

\begin{document}

\maketitle

\begin{abstract}
Based on the complete-lattice approach, a new Lagrangian duality theory for set-valued optimization problems is presented. In contrast to previous approaches, set-valued versions for the known scalar formulas involving infimum and supremum are obtained. In particular, a strong duality theorem, which includes the existence of the dual solution, is given under very weak assumptions: The ordering cone may have an empty interior or may not be pointed. "Saddle sets" replace the usual notion of saddle points for the Lagrangian, and this concept is proven to be sufficient to show the equivalence between  the existence of primal/dual solutions and strong duality on the one hand and the existence of a saddle set for the Lagrangian on the other hand.
\end{abstract} 

\pagebreak

%%%New Section
\section{Introduction}

This paper is concerned with a Lagrange duality theory for convex set-valued optimization with a set-valued constraint. In particular, a strong duality theorem is given, new concepts for saddle points in the set-valued framework are introduced and the relationship between primal/dual solutions and saddle points is established. The Lagrangian involves a set-valued analog of continuous linear functions generated by a pair of continuous linear functionals. This is a new construction and different from almost all known approaches since it is most common to use continuous linear operators as dual variables. Another important new feature is that our results involve (attained) infima/suprema with respect to set relations rather than only minimal/maximal elements. The theory given in \cite{Loehne11Book} is extended in several directions. Most notably, we do not assume that the ordering cone in the original image space has a non-empty interior. The main results also establish a one-to-one relationship between the set-valued problem and a corresponding problem for a family of scalarizations - despite the fact that we do not rely on so-called weak solutions.

Optimization problems with constraints given by set-valued functions ("correspondences", "relations", "multi-valued mappings") have been considered, for example, already in the late 70ies and early 80ies in a rather general framework. Compare the papers \cite{Borwein77}, \cite{Borwein81a}, \cite{Borwein81b} and \cite{Oettli81}, \cite{Oettli82}. The objective function was assumed to be extended real-valued or vector-valued (in \cite{Borwein81b}) or even set-valued (in \cite{Borwein81a}).

Duality theory for vector problems attracted attention to set-valued optimization problems, see \cite{TaninoSawaragi80}, \cite{Corley87}, \cite{Luc88}.  The paper \cite{TaninoSawaragi80} is particularly interesting since the authors defined a set-valued Lagrangian for a vector-valued problem: The Lagrangian is defined as a (negative) conjugate of a perturbation function with respect to the perturbation variable, and the supremum involved in the definition of the conjugate is replaced by a set of maximal points. A similar approach for a problem with a set-valued objective can be found in \cite{Postolica86}. In \cite{Tanino92}, a generalization of the supremum in $\R$ to a set-valued framework is used which later turned out to be a supremum in the sense of a complete lattice, compare \cite{LoehneTammer07}.

The introduction of so-called set relations (extensions of vector preorders to the power set of a preordered vector space) led to new solution concepts for set-valued optimization problems, see \cite{KuroiwaTanakaHa97}, \cite{Kuroiwa98-1,Kuroiwa03}, \cite{HernandezRodMarin07c}, \cite{Jahn04Book}. It has been pointed out by the authors of the present paper that set relations become even more valuable if they are used to construct appropriate image spaces for set-valued functions which carry the structure of complete lattices, see \cite{Hamel05Habil}, \cite{Loehne05Diss}, \cite{Hamel09}, \cite{Loehne11Book}. We obtain meaningful analogs to scalar formulas involving an infimum or a supremum. This feature is lacking in any vector optimization (duality) theory known to the authors and also in the set relation approach used in \cite{Ha05} , \cite{HernandezRodMarin07a}, \cite{HernandezRodMarin07b} where the aim is to find minimal elements in a set of sets with respect to a set relation. See also \cite{BotGradWanka09Book} for a discussion of Lagrange duality results in vector optimization with a set-valued objective.

The paper is organized as follows. In the next section, image spaces for set-valued functions are introduced, the basic problem is formulated and its Lagrangian defined. In the following sections, our solution concept for set-valued optimization problems is introduced, Fenchel conjugates are recalled and a scalarization concept is presented. The main results are contained in section \ref{SecDuality}. Saddle point results can be found in section \ref{SecSaddle}, and the appendix contains a definition of the concept "conlinear space".

%%%New Section
\section{Problem formulation, notation}

\label{SecProblem}

Let $Z$ be a topological linear space and $C \subseteq Z$ a convex cone with $0 \in C$.
We write $z_1\leq_C z_2$ for $z_2 - z_1 \in C$ with $z_1, z_2 \in Z$ which defines a
reflexive and transitive relation (a preorder). 
The relation $\leq_C$ on $Z$ can be extended to the powerset $\P(Z)$ of $Z$, the set of
all subsets of $Z$ including the empty set $\emptyset$, in two canonical ways (see \cite{Hamel09}
and the references therein). This gives rise to consider the following subsets of
$\P\of{Z}$:
\begin{align*}
\F\of{Z, C} & = \cb{A \in \P\of{Z} \;\vert\; A = \cl\of{A+C}} \\
 \label{EqDefQ} \G\of{Z, C} & = \cb{A \in \P\of{Z} \;\vert\; A = \cl\co\of{A+C}}.
\end{align*}
We shall abbreviate $\F\of{Z, C}$ and $\G\of{Z, C}$ to $\F\of{C}$ and $\G\of{C}$, respectively.

The Minkowski (elementwise) addition for non-empty subsets of $Z$ is extended to
$\P\of{Z}$ by
\[
\emptyset + A = A + \emptyset = \emptyset
\]
for $A \in \P\of{Z}$. Using this, we define an associative and commutative binary operation $\oplus \colon
\F\of{C} \times \F\of{C} \to \F\of{C}$ by
\begin{eqnarray}
\label{EqSumInP} A \oplus  B = \cl\of{A+B}
\end{eqnarray}
for $A, B \in \F\of{C}$. The element-wise multiplication of a set $A \subseteq Z$ with a
(non-negative) real number is extended by
\[
0 \cdot A = \cl C, \quad t \cdot \emptyset = \emptyset
\]
for all $A \in \F\of{C}$ and $t > 0$. In particular, $0 \cdot \emptyset = \cl C$ by
definition, and we will drop the $\cdot$ in most cases. The triple $\of{\F\of{C}, \oplus, \cdot}$ is a conlinear space with neutral element $\cl C$, and, obviously, $\of{\G\of{C}, \oplus, \cdot}$ is a conlinear subspace of it (see appendix for definitions).

On $\F\of{C}$ and $\G\of{C}$, $\supseteq$ is a partial order which is compatible with the
algebraic operations just introduced. Thus, $\of{\F\of{C}, \oplus, \cdot, \supseteq}$ and
$\of{\G\of{C}, \oplus, \cdot, \supseteq}$ are partially ordered, conlinear spaces in the sense
of \cite{Hamel05Habil}, \cite{Hamel09}. This property does not depend on properties of the cone $C$, in particular, $\leq_C$ does not need to
be a partial order. The use of $\supseteq$ can be motivated by considering equivalence
classes of an extension of the reflexive transitive relation on $Z$ generated by $C$ to
the power set of $Z$. Compare \cite{Hamel05Habil}, \cite{Hamel09} for more details.

Moreover, the pairs $\of{\F\of{C}, \supseteq}$, $\of{\G\of{C}, \supseteq}$ are complete lattices.
If $\emptyset \neq \mathcal A \subseteq \F\of{C}$, $\emptyset \neq \mathcal B \subseteq
\G\of{C}$, the following formulas hold true:
\[
 \inf_{\of{\F\of{C}, \supseteq}} \mathcal A = \cl\bigcup\limits_{A \in \mathcal A} A, \quad
 \sup_{\of{\F\of{C}, \supseteq}} \mathcal A = \bigcap\limits_{A \in \mathcal A} A
\]
and
\[
 \inf_{\of{\G\of{C}, \supseteq}} \mathcal B  = \cl\co\bigcup\limits_{B \in \mathcal B} B, \quad
 \sup_{\of{\G\of{C}, \supseteq}} \mathcal B = \bigcap\limits_{B \in \mathcal B} B.
\]
If $\mathcal A = \emptyset$ we set $\inf_{\of{\F\of{C}, \supseteq}} \mathcal A =
\emptyset$ and $\sup_{\of{\F\of{C}, \supseteq}} \mathcal A = Z$. This is in
accordance with the following monotonicity property: If $\mathcal A_1 \subseteq \mathcal
A _2$ then $\inf \mathcal A_1 \subseteq \inf \mathcal A_2$ and $\sup \mathcal A_1
\supseteq \sup \mathcal A_2$ in $\F\of{C}$. Likewise in $\G\of{C}$.

Let $X, Y$ be two locally convex spaces with topological duals $X^*, Y^*$, and $D \subseteq Y$ a convex
cone. The sets $\F\of{D} = \F\of{Y, D}$ and $\G\of{D} = \G\of{Y, D}$ are defined in the same way as $\F\of{C}$ and $\G\of{C}$.

Finally, let $f \colon X \to \F\of{C}$ and $g \colon X \to \F(D)$ be two functions. We
are interested in the problem
\[
\tag{P} \mbox{minimize} \quad f \quad \mbox{subject to} \quad  0 \in g\of{x}.
\]
The minimization is understood as looking for 
\[
p = \inf_{\of{\F\of{C}, \supseteq}}\cb{f\of{x} \;\vert\; x \in X, \; 0 \in g\of{x}} =
    \cl\bigcup\limits_{\cb{x \in X \;\vert\; 0 \in g\of{x}}} f\of{x},
\]
and a subset of minimizers in which this infimum is attained. A definition for the attainment of the infimum and a corresponding solution concept is given below. Likewise for functions mapping into $\G\of{C}$. This approach is different from most other approaches using a set relation such as \cite{Kuroiwa98-1, HernandezRodMarin07a}  which only focus on minimality notions.

To simplify our notation, we will drop the space and the order relations in expressions like $\inf_{\of{\F\of{C}, \supseteq}}\cb{f\of{x} \;\vert\; x \in X, \; 0 \in g\of{x}}$ in most cases below. By convention, an infimum is defined through the image space of the function: If $f$ maps into $\mathcal F\of{C}$, then the infimum is taken in $\of{\F\of{C}, \supseteq}$ and so on.

\begin{remark}
Since $g\of{x} \in \F\of{D}$ for each $x \in X$ we have
\[
0 \in g\of{x} \quad \Leftrightarrow \quad g\of{x} \cap -D \neq \emptyset.
\]
Already in \cite{Borwein77} it is shown that most of the common forms of constraints in the form of
inequalities and equations  (and many more) can be expressed in the above form. One may also compare chapter 5 of \cite{Luc88}.
\end{remark}

The positive dual (polar) cone of $C$ is $C^+ = \cb{z^* \in Z^* \colon \forall z \in C
\mid 0 \leq z^*\of{z}}$, and the negative dual is $C^- = -C^+$. We set $S\of{z^*}
= \cb{z \in Z \mid 0 \leq z^*\of{z}}$ for $z^* \in Z^*$.

Using the set-valued functions
\[
y \mapsto S_{\of{y^*, z^*}}\of{y} = \cb{z \in Z \mid y^*\of{y} \leq z^*\of{z}}
\]
we define the Lagrangian $l \colon X \times Y^* \times C^+\bs\cb{0} \to \F\of{C}$ of the
problem (P) by
\begin{equation}\label{DefLagrangian}
l\of{x, y^*, z^*} = f\of{x} \oplus \bigcup_{y \in g\of{x}}S_{\of{y^*, z^*}}\of{y}
    = f\of{x} \oplus \inf_{\of{\G\of{C}, \supseteq}}\cb{S_{\of{y^*, z^*}}\of{y} \colon y \in g\of{x}}.
\end{equation}
If $g\of{x} = \emptyset$, then the infimum over $y \in g\of{x} = \emptyset$ is the empty
set. The rules for the addition in $\of{\G\of{C}, \oplus}$ ($A \oplus \emptyset =
\emptyset$ for each $A \in \G\of{C}$) produce $l\of{x, y^*, z^*} = \emptyset$ in this
case.

Note that the functions $S_{\of{y^*, z^*}} \colon Y \to \mathcal P\of{Z}$ are positively homogeneous, additive with $S_{\of{y^*, z^*}}\of{0} = S\of{z^*}$, and they map into $\G\of{C}$ if and only if $z^* \in C^+$. If $z^* = 0$, then $S_{\of{y^*, z^*}}\of{y} \in \cb{Z, \emptyset}$, and otherwise $S_{\of{y^*, z^*}}\of{y}$ is a closed half space. In \cite{Hamel09}, it is shown that these functions provide a suitable substitute for continuous linear functions within the set-valued framework, and they have been applied in \cite{Loehne11Book} to obtain set-valued duality results for vector optimization problems (chapter 3, "duality of type II").

Under a mild condition, the primal problem can be reconstructed from the Lagrangian.

\begin{proposition}
\label{PropReconPrimal} {\bf (reconstruction of the primal)} If $f\of{x} \in \G\of{C}$,
$g\of{x} \in \G\of{D}$ and $f\of{x} \neq Z$ for each $x \in X$, then
\[
\sup_{\of{y^*, z^*} \in Y^* \times C^+\bs\cb{0}} l\of{x, y^*, z^*}
    = \bigcap_{\of{y^*, z^*} \in Y^* \times C^+\bs\cb{0}} l\of{x, y^*, z^*}
    = \left\{
    \begin{array}{ccc}
    f\of{x} & : & 0 \in g\of{x} \\
    \emptyset & : & 0 \not\in g\of{x}
    \end{array}
    \right.
\]
\end{proposition}

{\sc Proof.} First, assume $0 \in g\of{x}$. Then
\begin{align*}
\sup_{\of{y^*, z^*} \in Y^* \times C^+\bs\cb{0}} l\of{x, y^*, z^*}
    & \supseteq \sup_{\of{y^*, z^*} \in Y^* \times C^+\bs\cb{0}} f\of{x} \oplus S\of{z^*}\\
    & = \bigcap_{z^* \in C^+\bs\cb{0}} f\of{x} \oplus S\of{z^*} = f\of{x}
\end{align*}
since $f$ maps into $\G\of{C}$ and therefore $f\of{x}$ is the intersection of all closed
half spaces containing it. On the other hand,
\begin{align*}
\sup_{\of{y^*, z^*} \in Y^* \times C^+\bs\cb{0}} l\of{x, y^*, z^*}
    & \subseteq \bigcap_{z^* \in C^+\bs\cb{0}} l\of{x, 0, z^*} \\
    & = \bigcap_{z^* \in C^+\bs\cb{0}} f\of{x} \oplus S\of{z^*} = f\of{x}
\end{align*}
since $S_{\of{0, z^*}}\of{y} = S\of{z^*}$ and $f$ maps into $\G\of{C}$.

Next, assume $0 \not\in g\of{x}$. Since $g\of{x}$ is closed and convex, a separation
argument produces $y^* \in Y^*\bs\cb{0}$ such that
\[
\inf_{y \in g\of{x}} y^*\of{y} > 0.
\]
This implies $\lim\limits_{n \to \infty} \inf_{y \in g\of{x}} \of{ny^*}\of{y} = +\infty$,
hence for each $z^* \in C^+\bs\cb{0}$
\[
\bigcap_{n \in \N}\of{f\of{x} \oplus \cl\bigcup_{y \in g\of{x}}S_{\of{ny^*, z^*}}\of{y}}
    = \emptyset
\]
since
\[
\bigcup_{y \in g\of{x}}S_{\of{ny^*, z^*}}\of{y} \subseteq
    \cb{z \in Z \colon \inf_{y \in g\of{x}}\of{ny^*}\of{y} \leq z^*\of{z} }.
\]
This completes the proof. \pend

\begin{remark}
The condition $f\of{x} \in \G\of{C}$, $g\of{x} \in \G\of{D}$ for all $x \in X$ could be
understood as a well-posedness condition for the problem. If it is violated, then the
transition from the original problem to its Lagrangian form results in a loss of
information -- or one agrees to replace the original problem by
\[
\mbox{minimize} \quad \cl\co f\of{x} \quad \mbox{subject to} \quad  0 \in \cl\co g\of{x}.
\]
\end{remark}

%%%New section
\section{Solutions for set-valued optimization problems and canonical extensions}

In this section, we introduce a solution concept for complete lattice-valued problems as
well as saddle points for corresponding bi-variable functions. Let $f \colon X \to
\mathcal L$ be a function with values in a complete lattice $\of{\mathcal L, \leq}$.
Consider the two optimization problems
\begin{align}
 \tag{$\mathcal L_{min}$}\label{lvp}
  & \text{ minimize } f \colon X \to \mathcal L \; \text{ w.r.t. } \leq \text{ over } X ,\\
 \tag{$\mathcal L_{max}$}\label{lvp_max}
  & \text{ maximize } f \colon X \to \mathcal L \; \text{ w.r.t. } \leq \text{ over } X.
\end{align}
The sets of minimal and maximal elements of a subset $A \subseteq \mathcal L$ are defined
as usual by
\begin{align*}
\Min A & = \cb{z \in A \; \vert \; \of{y \in A \wedge y \leq z} \implies y=z}, \\
\Max A & = \cb{z \in A \; \vert \; \of{y \in A \wedge y \geq z} \implies y=z}.
\end{align*}

Moreover, for a set $M \subseteq X$ we denote the collection of all values of $f$ over
$M$ by
\[
f[M] = \cb{f\of{x} \vert \; x \in M} \subseteq  \mathcal L.
\]
A solution of ($\mathcal L_{min}$), for example, is expected to satisfy some minimality
condition. On the other hand, one also expects the infimum of $f$ to be attained at a
solution. In contrast to the set of extended real numbers, these two requirements do not
coincide in a general complete lattice. The next definition deals with minimal and
maximal function values.

\begin{definition}
\label{DefMiniMaximizer} An element $\bar x \in X$ is called a {\em minimizer} of $f$ if
$f\of{\bar x} \in \Min f\sqb{X}$, and it is called a {\em maximizer} of $f$ if $f\of{\bar
x} \in \Max f\sqb{X}$.
\end{definition}

Next, we deal with the attainment of the infimum/supremum.

\begin{definition}
A set $\bar X \subseteq X$ is called an {\em infimizer} of $f$ if
\begin{equation}\label{eq_inf_att}
 \inf_{x \in \bar X} f\of{x} = \inf_{x \in X}f\of{x},
\end{equation}
and in this case, we say that the infimum of $f$ is attained in $\bar X$. Likewise, a set
$\bar X \subseteq X$ is called a {\em supremizer} of $f$ if
\begin{equation}\label{eq_sup_att}
 \sup_{x \in \bar X} f\of{x} = \sup_{x \in X}f\of{x},
\end{equation}
and in this case, we say that the supremum of $f$ is attained in $\bar X$.
\end{definition}

One cannot expect that an infimizer or supremizer is a singleton. Indeed, if, for
example, we consider a function $f$ into the complete lattice $\of{\F\of{C}, \supseteq}$,
the infimum is given by the closure of the union of the function values which is not a
function value itself in general. Therefore, we agree that a solution of ($\mathcal
L_{min}$) should be a subset of $X$ rather than a single element. This is the point of
the following definition which is essentially due to \cite{HeydeLoehne10}.

\begin{definition}\label{def_sol}
An infimizer $\bar X \subseteq X$ of $f$ is called a {\em solution} to \eqref{lvp} if
$f\sqb{\bar X} \subseteq \Min f\sqb{X}$, and it is called a {\em full solution} if $
f\sqb{\bar X} = \Min f\sqb{X}$. Similarly, a supremizer $\bar X \subseteq X$ of $f$ is
called a {\em solution} to \eqref{lvp_max} if $f\sqb{\bar X} \subseteq \Max f\sqb{X}$,
and it is called a {\em full solution} if $ f\sqb{\bar X} = \Max f\sqb{X}$.
\end{definition}

Clearly, a solution to \eqref{lvp} is an infimizer which consists of only minimizers.
Note that the earliest approach to set-valued optimization problems consisted in looking
for minimal points of the union of all function values, see \cite{Corley87}, \cite{Luc88}, for example. Thus, one
could understand this as looking for a minimal element of the infimum. On the other hand,
the set relation approach as presented for instance in \cite{Kuroiwa98-1, HernandezRodMarin07c} aims at finding minimizers with
respect to a set relation without caring for the infimum.

\begin{remark}
Note that a solution in the sense of Definition \ref{def_sol} above is called ''mild solution'' in \cite{HeydeLoehne10} and \cite{Loehne11Book} and a full solution is just called ''solution'' in these references. We decided to change the notation since, in particular, in linear vector optimization problems one wants to have that a solution consists of finitely many elements.
\end{remark} 

\begin{example}
Let the objective function $f \colon \R \to \G\of{\R^2, \R^2_+}$ be defined by 
\[
f\of{x} = \cb{z \in \R^2_+ \vert \; z_1 \geq 3 + 2x + r,\; z_2 \geq 3 + 2x - r,\; r \in \sqb{-x^2,x^2}}
\]
for $x \in \R_+$ and $f\of{x} = \emptyset$ for $x \not\in \R_+$. The set $\bar X = \cb{0}
\cup (2,3]$ is a full solution to \eqref{lvp}, and the set $\bar X = \cb{0,3}$ is a
solution to \eqref{lvp}. This example shows that a full solution may entail a far larger
set $\bar X$ than being required for the attainment of the infimum.
\end{example}

The above solution concept can also be expressed in terms of the {\em canonical
extension} of the function $f$. We introduce this concept below since we will make use of it
when discussing saddle points.

\begin{definition}
\label{DefCanonicalExt} Let $f \colon X \to \mathcal L$ be a function with values in the
complete lattice $\of{\mathcal L, \leq}$. The function $\hat{F} \colon
\mathcal P\of{X} \to \mathcal L$ defined by
\[
\hat{F}\of{M} = \inf_{x \in M} f\of{x}
\]
is called the {\em (canonical) inf-extension} of $f$, whereas the function
$\check{F} \colon \mathcal P\of{X} \to \mathcal L$ defined by
\[
\check{F}\of{M} = \sup_{x \in M} f\of{x}
\]
is called the {\em (canonical) sup-extension} of $f$.
\end{definition}

Clearly, a set $\bar X \subseteq X$ is an infimizer of $f$ if and only if the
inf-extension $\hat{F}$ of $f$ attains its infimum in $\bar X$, that is
\[
\hat{F}\of{\bar X} = \inf_{M \in \mathcal
P\of{X}}\hat{F}\of{M},
\]
and likewise for the supremum of $\check{F}$. Therefore, the inf-extension will
play a role for the minimization problem ($\mathcal L_{min}$) and the sup-extension for
the maximization problem ($\mathcal L_{max}$).

We turn to the definition of saddle points for a function $l \colon X \times V \to
\mathcal L$. The following definition extends the concept of canonical extensions to such
functions. Note that there are two such extensions which seems natural since one can
consider "$\inf\sup$" as well as "$\sup\inf$".

\begin{definition}
\label{DefCanExtension2D} Let $l \colon X \times V \to \mathcal L$ be a function. The
{\em lower canonical extension} of $l$ is the function $\hat{L} \colon
\mathcal P\of{X} \times \mathcal P\of{V} \to \mathcal L$ defined by
\[
\hat{L}\of{U, W} = \sup_{v \in W} \inf_{x \in U} l\of{x, v}.
\]
The {\em upper canonical extension} of $l$ is the function $\check{L} \colon
\mathcal P\of{X} \times \mathcal P\of{V} \to \mathcal L$ defined by
\[
\check{L}\of{U, W} = \inf_{x \in U}\sup_{v \in W} l\of{x, v}.
\]
%In both cases, $U \subseteq X$, $V \subseteq V$.
\end{definition}

The very definitions of $\hat{L}$, $\check{L}$ immediately produce
\begin{equation}
\label{EqLowerLessUpper}
 \forall U \in \mathcal P\of{X}, \; \forall W \in \mathcal P\of{V}
 \colon \hat{L}\of{U, W} \leq \check{L}\of{U, W}.
\end{equation}

\section{Fenchel conjugates for set-valued functions}

Given $x^* \in X^*$ and $z^* \in C^+$,  the (conlinear) function $S_{\of{x^*, z^*}} \colon X
\to \G(C)$ is defined as its corresponding version on $Y$: $S_{\of{x^*, z^*}}\of{x} = \cb{z \in Z \mid x^*\of{x} \leq z^*\of{z}}$
for $x \in X$.

The (negative) Fenchel conjugate of a function $f \colon X \to \mathcal P\of{Z}$ is the
function $-f^* \colon X^* \times \of{C^+\bs\cb{0}} \to \G(C)$ defined by
\begin{equation} \label{LFC} -f^*\of{x^*, z^*} =
    \cl\bigcup_{x \in X} \sqb{f\of{x} + S_{\of{x^*, z^*}}\of{-x}}.
\end{equation}
See \cite{Hamel09} for more details about and a motivation for this definition. We
immediately conclude the Young--Fenchel inequality for set-valued functions
\begin{equation}
\label{EqYFI} \forall z^* \in C^+\bs\cb{0}, \; \forall x \in X, \;  \forall x^*\in X^*  \colon
    -f^*\of{x^*, z^*} \supseteq f\of{x} + S_{\of{x^*, z^*}}\of{-x}.
\end{equation}
Note that the values of $-f^*$ are closed half spaces or $\emptyset$ or $Z$. The main
result about Fenchel conjugates is the Fenchel--Moreau theorem which establishes
conditions under which the biconjugate of a function $f$ coincides with $f$. The
biconjugate $f^{**} \colon X \to \G(C)$ of $f \colon X \to \G(C)$ is
given by
\[
f^{**}\of{x} = \bigcap_{x^* \in X^*, \; z^* \in \of{C^+\bs\cb{0}}}\sqb{-f^*\of{x^*, z^*}
+ S_{\of{x^*, z^*}}\of{x}}.
\]
The graph of a function $f \colon X \to \mathcal P\of{Z}$ is the set
\[
\gr f = \cb{\of{x,z} \in X \times Z \mid z \in f\of{x}} \subseteq X \times Z.
\]
A function $f \colon X \to \mathcal P\of{Z, C}$ is called convex if its graph is a convex set. This is equivalent to 
\[
\forall x_1, x_2 \in X, \forall t \in \of{0,1} \colon f\of{tx_1 + \of{1-t}x_2} \supseteq tf\of{x_1} + \of{1-t}f\of{x_2}.
\]
We define the closure $\cl f$ of $f$ by
\[
z \in \of{\cl f}\of{x} \quad \Leftrightarrow \quad \of{x, z} \in \cl\of{\gr f}.
\]
With this definition, we can state the Fenchel--Moreau theorem in the following form.

\begin{theorem}
\label{ThmFM} Let $f \colon X \to \G(C)$ be a convex function such that there is
$x_0 \in X$ with $f\of{x_0} \neq \emptyset$ and $\of{\cl f}\of{x} \neq Z$ for all $x \in
X$. Then $\cl f = f^{**}$.
\end{theorem}

{\sc Proof.} See \cite{Hamel09}, Theorem 3. \pend

\section{Linear scalarization and the scalarized problems}

Again, under some well-posedness assumptions, a function $f \colon X \to \mathcal
P\of{Z}$ and the set-valued problem (P) can equivalently be described by families of
extended real-valued functions and scalar problems, respectively. The scalarization
approach has been used as the main tool for obtaining set-valued duality results of
Fenchel type along with a complete calculus for set-valued convex functions in
\cite{Schrage09Diss}, see also \cite{HamelSchrage12} and \cite{Schrage10R}. The
set-valued approach to vector optimization summarized in \cite{Loehne11Book} also relies
on a scalarization procedure. Here, we additionally show the equipollence of the
set-valued results with corresponding results for families of scalarizations.

\subsection{Scalarization of $\G(C)$-valued functions}

Let a function $f \colon X \to \mathcal P\of{Z}$ and $z^* \in C^+$ be given. Define an extended real-valued
function $\vp_{f, z^*} \colon X \to \OLR = \R\cup\cb{\pm\infty}$ by
\begin{equation}
\label{EqScalarize} \vp_{f, z^*}\of{x} = \inf_{z \in f\of{x}}z^*\of{z}.
\end{equation}
The set-valued function $f \colon X \to \P(C)$ is convex if and only if $\vp_{f,
z^*}\colon X \to \OLR$ is convex for each $z^* \in C^+$. Of course, $\vp_{f, z^*}\of{x}$
is nothing else than the value of negative support function of the set $f\of{x}$ at
$-z^*$, but we shall emphasize the dependence on $x$ rather than on $z^*$. The support
function interpretation of $\vp_{f, z^*}$ immediately gives the following formula. If $f
\colon X \to \F\of{C}$ is convex, then $f\of{x} \in \G(C)$ for all $x \in X$ and
\begin{equation}
\label{EqReScalarize} \forall x \in X \colon f\of{x} =
    \bigcap_{z^* \in C^+\bs\cb{0}} \cb{z \in Z \mid \vp_{f, z^*}\of{x} \leq z^*\of{z}}.
\end{equation}
This is a consequence of the separation theorem since the closed convex set $f\of{x}$ is
the intersection of all closed half spaces containing it, and only half spaces with
normals $z^* \in C^+\bs\cb{0}$ need to be considered since $f\of{x} \in \G(C)$.

Likewise, the extended real-valued functions $\vp_{g, y^*} \colon X \to \OLR$ with $y^*
\in D^+\bs\cb{0}$ are defined for the set-valued function $g \colon X \to \P(D)$.

The next result shows that scalarization and conjugation commute.

\begin{lemma}
\label{LemScalarConjugate} Let $f \colon X \to \mathcal P\of{Z}$ be a function and $x^* \in X^*$, $z^* \in C^+\bs\cb{0}$. Then
\begin{align*}
-\vp^*_{f, z^*}\of{x^*} & = \inf_{z \in -f^*\of{x^*, z^*}} z^*\of{z}, \\
\of{-f^*}\of{x^*, z^*} & = \cb{z \in Z \mid -\vp^*_{f, z^*}\of{x^*} \leq z^*\of{z}}
\end{align*}
where $\vp^*_{f, z^*} = \of{\vp_{f, z^*}}^* \colon X \to \OLR$ is the classical (scalar) Fenchel
conjugate of $\vp_{f, z^*}$. 
\end{lemma}

{\sc Proof.} The proof of the first equation is straightforward from the definitions of $-\vp_{f, z^*}$ and (scalar) Fenchel conjugates. The second is obviously true if $f \equiv \emptyset$ $\Leftrightarrow$ $-f^* \equiv \emptyset$. If this is not the case, take $z \in  f\of{x} + S_{\of{x^*, z^*}}\of{-x}$. Then there are $z_1 \in f\of{x}$ and $z_2 \in S_{\of{x^*, z^*}}\of{-x}$ such that $z = z_1 + z_2$. Using the
definition of $S_{\of{x^*, z^*}}$ we obtain $-x^*\of{x} \leq z^*\of{z_2}$. From the
definition of $\vp_{f, z^*}$ we get $\vp_{f, z^*}\of{x} \leq z^*\of{z_1}$. Hence
\[
\vp_{f, z^*}\of{x} -x^*\of{x} \leq z^*\of{z_1 + z_2} = z^*\of{z}
\]
and therefore
\[
-\vp^*_{f, z^*}\of{x^*} = \inf_{x \in X}\sqb{\vp_{f, z^*}\of{x} -x^*\of{x}} \leq z^*\of{z}.
\]
This shows $\bigcup_{x \in X}\sqb{f\of{x} + S_{\of{x^*, z^*}}\of{-x}} \subseteq \cb{z \in Z \mid -\vp^*_{f, z^*}\of{x^*} \leq z^*\of{z}}$. Since
$z^*$ is a continuous function, this implies $-f^*\of{x^*, z^*} \subseteq \cb{z \in
Z \mid -\vp^*_{f, z^*}\of{x^*} \leq z^*\of{z}}$.

Conversely, take $z_0 \in Z$ satisfying $-\vp^*_{f, z^*}\of{x^*} \leq z^*\of{z_0}$. If $z_0
\not\in -f^*\of{x^*, z^*}$ there would be a $z^*_0 \in Z^*\bs\cb{0}$ such that
\[
z^*_0\of{z_0} < \inf_{z \in -f^*\of{x^*, z^*}} z^*_0\of{z}.
\]
Since $z^*_0 \neq 0$ and $-f^*\of{x^*, z^*}$ is a closed half space with normal $z^*$ we can assume, without loss of generality, $z^*_0 = z^*$ and obtain the contradiction $z^*\of{z_0} < -\vp^*_{f, z^*}\of{x^*} \leq z^*\of{z_0}$. \pend

\subsection{Scalarization of the optimization problem}
\label{SubSecScalarPrimal}

Together with (P) we consider the family of scalar problems
\[
\tag{PS}  \mbox{minimize} \quad \vp_{f, z^*} \quad \mbox{subject to} \quad  0 \in g\of{x}
\]
for $z^* \in C^+\bs\cb{0}$, that is we want to find the numbers
\[
p_{z^*} = \inf\cb{\vp_{f, z^*}\of{x} |\; 0 \in g\of{x}}, \; z^* \in C^+\bs\cb{0}.
\]
Together with the Lagrange function for (P) as introduced in \eqref{DefLagrangian} we consider the Lagrange function of the scalarized problem
\[
\lambda_{z^*}\of{x, y^*} = \vp_{f, z^*}\of{x} + \inf_{y \in g\of{x}}y^*\of{y} =
    \vp_{f, z^*}\of{x} + \vp_{g, y^*}\of{x}.
\]
The set of functions
\[
\cb{\lambda_{z^*}\of{\cdot, \cdot} |\; X \times Y^* \to \OLR \colon z^* \in
C^+\bs\cb{0}}
\]
is called the (scalarized) Lagrangian family of problem (P). One main message of this
paper is that every Lagrange duality result about (P) can equivalently be expressed in
terms of its Lagrangian family.

\begin{proposition}
\label{PropLagrangeScalarized} {\bf (Lagrange functions of the scalarized problems)} We
have
\[
\lambda_{z^*}\of{x, y^*} = \vp_{l, z^*}\of{x, y^*}
\]
for all $x\in X$, $y^* \in Y^*$, $z^* \in C^+$.
\end{proposition}

{\sc Proof.} We have
\begin{align*}
\vp_{l, z^*}\of{x, y^*} & = \inf_{z \in l\of{x, y^*, z^*}} z^*\of{z} \\
    & = \inf_{z \in f\of{x}}z^*\of{z} + \inf_{z \in \cl\bigcup_{y \in g\of{x}}S_{\of{y^*,
    z^*}}\of{y}} z^*\of{z} \\
    & = \vp_{f, z^*}\of{x} +
    \inf\cb{z^*\of{z} \mid z \in \bigcup_{y \in g\of{x}}S_{\of{y^*, z^*}}\of{y}}
\end{align*}
by calculus rules for support functions and continuity of $z^*$. Since
\[
\vp_{S_{\of{y^*, z^*}}, z^*}\of{y} = y^*\of{y}
\]
by definition of $S_{\of{y^*, z^*}}$ we may conclude
\[
\vp_{l, z^*}\of{x, y^*}  = \vp_{f, z^*}\of{x} + \inf\cb{y^*\of{y} \mid y \in g\of{x}},
\]
and the right hand side is just $\lambda_{z^*}\of{x, y^*}$. \pend

\medskip The previous proposition tells us that the two operations "scalarization" and
"transition to the Lagrangian" commute.

\begin{proposition}
\label{PropInfScalarization} {\bf (scalarized problems)} Let $f\colon X \to \P(C)$
and $g\colon X \to \P(D)$. The following statements are equivalent for $z^* \in
C^+$:
\\
(a) $\cl\co\bigcup\cb{f\of{x} \oplus S\of{z^*} \mid 0 \in g\of{x}} \neq Z$
\\
(b) $-\infty < \inf\cb{\vp_{f, z^*}\of{x} \mid 0 \in g\of{x}}$

\end{proposition}

{\sc Proof.} Strainghtforward. \pend

%%%Newsection
\section{Duality}
\label{SecDuality}

%%%Newsubsection
\subsection{Construction of the dual problem and weak duality}
\label{SubsecDual}

As in scalar optimization, the objective of the dual problem is the infimum of the
Lagrangian with respect to the primal variable, that is the function $h \colon Y^* \times
C^+\bs\cb{0} \to \G\of{C}$ defined by
\[
h\of{y^*, z^*} = \inf_{x\in X}l\of{x, y^*, z^*} = \cl\bigcup_{x \in X} l\of{x, y^*, z^*}.
\]
Since the values of $l$ are closed half spaces, the convex hull can be dropped in the
infimum. The dual problem,
\[
\tag{D} \mbox{maximize} \quad h \quad \mbox{subject to} \quad  y^* \in Y^*, \;z^* \in
 C^+\bs\cb{0},
\]
thus consists in finding
\[
d = \sup_{y^* \in Y^*, \, z^* \in C^+\bs\cb{0}} h\of{y^*, z^*} =
    \bigcap _{y^* \in Y^*, \, z^* \in C^+\bs\cb{0}} h\of{y^*, z^*}
\]
and corresponding (full) solutions.

The objectives of the scalarized dual problems (see Section \ref{SubSecScalarPrimal}) are
\[
\inf_{x \in X}\lambda_{z^*}\of{x, y^*} = \inf_{x\in X}\sqb{\vp_{f, z^*}\of{x} + \vp_{g,
y^*}\of{x}}.
\]
The scalarized dual problems become
\[
\tag{DS} \mbox{maximize} \quad \inf_{x \in X}\lambda_{z^*}\of{x, \cdot} \quad
\mbox{subject to} \quad  y^* \in Y^*,
\]
that is find
\[
d_{z^*} = \sup\cb{\inf_{x \in X}\lambda_{z^*}\of{x, y^*} \mid y^* \in Y^*}, \; z^* \in
C^+\bs\cb{0}
\]
and corresponding solutions.

\begin{proposition}
Scalarization and dualization commute, that is
\[
\vp_{h, z^*}\of{y^*} = \inf_{x\in X}\lambda_{z^*}\of{x, y^*}
\]
for all $y^* \in Y^*$ and $z^* \in C^+\bs\cb{0}$.
\end{proposition}

{\sc Proof.} Straightforward. \pend

\medskip The following weak duality result is immediate.

\begin{proposition}
\label{PropWeakDuality} {\bf (weak duality)} The following equivalent statements are
true:
\\
(a) For all $x \in X$ satisfying $0 \in g\of{x}$, for all $\of{y^*, z^*} \in Y^* \times
C^+\bs\cb{0}$,
\[
h\of{y^*, z^*} \supseteq f\of{x} \oplus S\of{z^*}.
\]
(b) For all $x \in X$ satisfying $0 \in g\of{x}$, for all $\of{y^*, z^*} \in Y^* \times
C^+\bs\cb{0}$,
\[
\vp_{h, z^*}\of{y^*} = \inf_{z \in h\of{y^*, z^*}}z^*\of{z} \leq \inf_{z \in
f\of{x}}z^*\of{z} = \vp_{f, z^*}\of{x}.
\]
\end{proposition}

{\sc Proof.} (a) is true since for $\of{y^*, z^*} \in Y^* \times C^+\bs\cb{0}$ and $x \in
X$ satisfying $0 \in g\of{x}$ we have
\begin{align*}
h\of{y^*, z^*} & \supseteq f\of{x} \oplus
        \cl\bigcup_{y \in g\of{x}}S_{\of{y^*, z^*}}\of{y} \\
    & \supseteq f\of{x} \oplus S_{\of{y^*, z^*}}\of{0} = f\of{x} \oplus S\of{z^*}.
\end{align*}
The implication (a) $\Rightarrow$ (b) is obvious. The converse is also true since
$h\of{y^*, z^*}$ and $f\of{x} \oplus G\of{z^*}$ are half spaces with normal $z^* \neq 0$
and $\inf_{z \in f\of{x}}z^*\of{z} = \inf_{z \in f\of{x}}z^*\of{z} \oplus S\of{z^*}$.
\pend

\subsection{Strong duality}

The value function of the basic optimization problem is the function  $v \colon Y \to
\mathcal \G\of{C}$ defined by
\[
v\of{y} = \inf\cb{f\of{x} \mid x \in X,\; y \in g\of{x}}.
\]
Using again the convention that the infimum over an empty set is $\emptyset$, we obtain
$v\of{y} = \emptyset$ whenever $\cb{x \in X \mid y \in g\of{x}} = \emptyset$. Clearly,
\[
p = v\of{0} = \inf\cb{f\of{x} \mid 0 \in g\of{x}}
\]
is the optimal value of the original problem.

The optimal values of the scalarized problems are defined by
\[
v_{z^*}\of{y} = \inf\cb{\vp_{f, z^*}\of{x} \mid y \in g\of{x}}, \; z^* \in
C^+\bs\cb{0}.
\]
Clearly,
\[
p_{z^*} = v_{z^*}\of{0} = \inf\cb{\vp_{f, z^*}\of{x} \mid 0 \in g\of{x}}, \; z^* \in
C^+\bs\cb{0},
\]
are the optimal values of the scalarized problems.

\begin{proposition}\label{PropScalarizedValue}
Scalarization and transition to value functions commute, that is
\[
\forall y \in Y \colon \vp_{v, z^*}\of{y} = \inf\cb{\vp_{f, z^*}\of{x} \mid y \in
g\of{x}} = v_{z^*}\of{y}
\]
with $\vp_{v, z^*} \colon Y \to \OLR$ being the scalarization of the value function.
\end{proposition}

{\sc Proof.} Obvious by construction. \pend

\begin{proposition}
\label{PropDualOptimalValue} (a) The (negative) Fenchel conjugate of the value function
is
\[
-v^*\of{y^*, z^*} = \inf\cb{l\of{x, -y^*, z^*} \mid x \in X},\]
and the negative Fenchel conjugate of the value function of the scalarized problem is
\[
-\of{v_{z^*}}^*\of{y^*} = \inf\cb{\lambda_{z^*}\of{x, -y^*} \mid x \in X}.\]

(b) $v^{**}\of{0} = d$, and $\of{v_{z^*}}^{**}\of{0} = d_{z^*}$ for all $z^* \in
C^+\bs\cb{0}$.
\end{proposition}

{\sc Proof.} (a) We have
\begin{align*}
-v^*\of{y^*, z^*} & = \cl\bigcup_{y \in Y} \sqb{v\of{y} + S_{\of{y^*, z^*}}\of{-y}} \\
    & = \cl \bigcup_{x \in X,\, y \in g\of{x}}
        \sqb{f\of{x} + S_{\of{y^*, z^*}}\of{-y}} \\
    & = \cl\bigcup_{x \in X} l\of{x, -y^*, z^*}\\
    & = \inf_{x \in X} l(x,-y^*,z^*).
\end{align*}

(b) Using the definition of $h$, part (a) and the properties of $-v^*$ we obtain
\begin{align*}
\bigcap_{y^*, \, z^*} h\of{y^*, z^*} &= \bigcap_{y^*, \, z^*}\cl\bigcup_{x \in X} l\of{x, y^*, z^*} \\
    &= \bigcap_{y^*, \, z^*}\sqb{-v^*\of{-y^*, z^*} + S_{\of{-y^*, z^*}}\of{0}}
    = v^{**}\of{0}.
\end{align*}
This completes the proof. \pend

\begin{proposition}
\label{PropZeroDualityGap} Let $f$ and $g$ be convex and $v\of{0} \neq Z$, $v\of{0} \neq
\emptyset$. Then, there is $z^* \in C^+\bs\cb{0}$ such that $p \oplus S\of{z^*} = v\of{0} \oplus S\of{z^*} \neq Z$, and the  following statements are equivalent:
\\
(a) $p = d$,
\\
(b) $v$ is closed at $0 \in Y$, that is $\of{\cl v}\of{0} = v\of{0}$,
\\
(c) $p_{z^*}=d_{z^*}$ for each $z^* \in C^+\bs\cb{0}$ with $p \oplus S\of{z^*} \neq Z$,
\\
(d) For each $z^* \in C^+\bs\cb{0}$ with $p \oplus S\of{z^*} \neq Z$ the function $y \mapsto \vp_{v, z^*}\of{y}$
is closed at $0 \in Y$, i.e. $\cl\of{\vp_{v, z^*}}\of{0} = \vp_{v, z^*}\of{0}$.
\end{proposition}

{\sc Proof.} 
Since $f$ and $g$ are convex, $v$ is convex as well. Moreover, $v\of{0} = p$ is
the optimal value of the original problem and a closed convex set. If $v\of{0} \not\in \cb{Z, \emptyset}$, a separation argument produces $z^* \in C^+\bs\cb{0}$ such that $v\of{0} \oplus S\of{z^*} \neq Z$.

The equivalence of (a) and (b) follows from (b) of Proposition
\ref{PropDualOptimalValue} and the set-valued Fenchel-Moreau theorem (Theorem \ref{ThmFM}).

The equivalence of (c) and (d) follows likewise, but it is also the known scalar result. 

It remains to prove the equivalence of (a) and (c). The implication (a) $\Rightarrow$ (c)
follows from the definition of $p$, $d$ and the scalarization functions. The converse is
a consequence of the fact that both $p$ and $d$ are closed convex sets, and one can apply
formula \eqref{EqReScalarize} with $f\of{x}$ replaced by $p$ and $d$, respectively. \pend

\medskip The remaining part of this section is devoted to sufficient conditions for
strong duality, or $v\of{0}= v^{**}\of{0}$. The following condition is called the {\bf
Slater condition} for problem (P):
\[
\exists \bar{x} \in \dom f \colon g\of{\bar{x}} \cap \Int\of{-D} \neq \emptyset.
\]
The implicit assumption is $\Int D \neq \emptyset$. 

\begin{theorem}\label{ThmStrongDuality}
{\bf (strong duality)} Let $f$ and $g$ be convex. If the Slater condition for problem (P) is satisfied, then strong
duality holds, that is the following two equivalent conditions are satisfied:
\begin{align}
p & = \inf\cb{f\of{x} \mid 0 \in g\of{x}}
    = \sup_{y^* \in Y^*, \, z^* \in C^+\bs\cb{0}} h\of{y^*, z^*} = d \\
\forall z^* \in C^+\bs\cb{0} \colon p_{z^*} & = \inf\cb{\vp_{f, z^*}\of{x} \mid 0 \in g\of{x}} = \sup_{y^* \in Y^*} \vp_{h, z^*}\of{y^*}  = d_{z^*}.
\end{align}
If the infimum of the primal problem is not $Z$, then the following two equivalent conditions are satisfied
\begin{align}
p \oplus S\of{z^*}\neq Z & \Rightarrow \; \exists y^* \in Y^* \colon p \oplus S\of{z^*} = h\of{y^*, z^*} \\
-\infty < p_{z^*}  & \Rightarrow \; \exists y^* \in Y^* \colon p_{z^*} = \vp_{h, z^*}\of{y^*}. 
\end{align}
Moreover, the set
\begin{multline*}
\Delta = \cb{\of{y^*, z^*} \in Y^* \times C^+\bs\cb{0} \mid Z \neq p \oplus S\of{z^*} = h\of{y^*, z^*}} \\ 
    = \cb{\of{y^*, z^*} \in Y^* \times C^+\bs\cb{0} \mid -\infty < p_{z^*}=  \vp_{h, z^*}\of{y^*}} 
\end{multline*}
is non-empty and a full solution of the dual problem.
\end{theorem}

The theorem can be summarized as "strong duality holds and is equivalent to strong duality for every meaningful scalarized problem".

\medskip {\sc Proof.} If $p=Z$, weak duality implies $d=p$. Thus let $p\neq Z$.  In view of Proposition
\ref{PropZeroDualityGap} we have to show that $v\of{0} =
\of{\cl v}\of{0}$. By Slater's condition there is a
neighborhood $V$ of $0 \in Y$ such that $y \in g\of{\bar{x}} \neq \emptyset$ for each $y
\in V$. Hence
\[
v\of{y} = \cl\bigcup\cb{f\of{x} \mid y \in g\of{x}} \supseteq
    f\of{\bar{x}} \neq \emptyset
\]
for all $y \in V$. Take $z_0 \in
\of{\cl v}\of{0}$, that is $\of{0, z_0} \in \cl\of{\gr v}$. Then, for each neighborhood
$W$ of $0 \in Y$ and each neighborhood $U$ of $0 \in Z$ we have
\[
\sqb{\of{0, z_0} + \of{W \times U}} \cap \gr v \neq \emptyset.
\]
Fix a neighborhood $U$ of $0 \in Z$ and $0 < t < 1$. Then $\sqb{\of{0, z_0} + \of{tV
\times U}} \cap \gr v \neq \emptyset$. Take $y_0 \in V$ and $z \in U$ such that $\of{0,
z_0} + \of{ty_0, z} \in \gr v$, and take $\bar z \in f\of{\bar x}$. Then
\[
\of{0, z_0} + \of{ty_0, z} \in \gr v \; \wedge \; \of{0, 0} + \of{-ty_0, \bar z} \in \gr
v.
\]
The second inclusion is a consequence of $f\of{\bar x} \subseteq v\of{y}$ for all $y \in
V$ since we can assume that $V$ is balanced. Taking the convex combination of the two
points in $\gr v$ with $s = \frac{1}{1+t}$ and $1-s = \frac{t}{1+t}$ we obtain
\[
\frac{1}{1+t}\of{z_0 + z} + \frac{t}{1+t}\bar z \in v\of{0}
\]
for all $t \in \sqb{0, 1}$. Since $v\of{0}$ is a closed set by definition, we let $t \to
0$ and get $z_0 + z \in v\of{0}$. Hence, for each neighborhood $U$ of $0 \in Z$ we have
$\of{z_0 + U} \cap v\of{0} \neq \emptyset$ which implies $z_0 \in v\of{0}$ since the
latter set is closed. This proves that $v\of{0} = \of{\cl v}\of{0}$ and hence $p=d$.

Next, since $p \in \mathcal \G\of{C}\setminus Z$ there is $z^* \in C^+\bs\cb{0}$ such that $p \oplus
G\of{z^*} \neq Z$, and we have $p \oplus S\of{z^*} = \cb{z \in Z \mid v_{z^*}\of{0}
\leq z^*\of{z}}$. Proposition \ref{PropScalarizedValue} implies $v_{z^*}\of{0} = \vp_{v,
z^*}\of{0}$, hence
\[
z \in p \oplus S\of{z^*} \quad \Leftrightarrow \quad \vp_{v, z^*}\of{0} \leq z^*\of{z}.
\]
Proposition \ref{PropInfScalarization} yields that the convex function $y \to
v_{z^*}\of{y} = \vp_{v, z^*}\of{y}$ has a finite value at $0 \in Y$. The Slater condition
ensures that this function is continuous at $0 \in Y$, hence there is $y^* \in Y^*$ such
that
\[
-y^* \in \partial\vp_{v, z^*}\of{0}
\]
which in turn is equivalent to
\[
\of{\vp_{v, z^*}}^*\of{-y^*} + \vp_{v, z^*}\of{0} = -y^*\of{0} = 0
\]
(Young-Fenchel equality for elements of the subdifferential). Finally, observe that
\[
z \in -v^*\of{y^*, z^*}  \quad \Leftrightarrow \quad
    -\of{\vp_{v, z^*}}^*\of{y^*} \leq z^*\of{z}.
\]
Altogether, we obtain
\begin{align*}
z \in p \oplus G\of{z^*} &
    \quad \Leftrightarrow \quad \vp_{v, z^*}\of{0} \leq z^*\of{z} \\
    & \quad \Leftrightarrow \quad -\of{\vp_{v, z^*}}^*\of{-y^*} \leq z^*\of{z}\\
    & \quad \Leftrightarrow \quad z \in -v^*\of{-y^*, z^*} \\
    & \quad \Leftrightarrow \quad z \in h\of{y^*, z^*}
\end{align*}
where the last equivalence is a consequence of Proposition \ref{PropDualOptimalValue},
(a). This proves that the set $\Delta$ is non-empty.

Finally, in view of Definition \ref{def_sol} we have to establish the following facts:
\\
(i) $d = \bigcap\cb{h\of{y^*, z^*} \mid \of{y^*, z^*} \in \Delta}$ (attainment of the supremum).
\\
(ii) Each $h\of{\bar{y}^*, \bar{z}^*}$ with $\of{\bar{y}^*, \bar{z}^*} \in \Delta$ is a
maximal element of
\[
H = \cb{h\of{y^*, z^*} \neq Z \mid \of{y^*, z^*} \in Y^* \times C^+\bs\cb{0}}
\]
with respect to $\supseteq$ ($\Delta$ is a solution since each of its elements is a
maximizer).
\\
(iii) If $h\of{\bar{y}^*, \bar{z}^*}$ is a maximal element of $H$, then $\of{\bar{y}^*,
\bar{z}^*} \in \Delta$ ($\Delta$ is a full solution since all maximizers of the dual
problem are elements of $\Delta$).

(i) is a direct consequence of strong duality.

(ii) $h\of{\bar{y}^*, \bar{z}^*}$ is maximal in $H$ if and only if
\[
\of{y^*, z^*} \in Y^* \times C^+\bs\cb{0}, \; h\of{\bar{y}^*, \bar{z}^*} \supseteq
h\of{y^*, z^*} \quad \Rightarrow \quad h\of{\bar{y}^*, \bar{z}^*} = h\of{y^*, z^*}.
\]
If $\of{\bar{y}^*, \bar{z}^*} \in \Delta$ and $h\of{\bar{y}^*, \bar{z}^*} \supseteq
h\of{y^*, z^*}$, then by weak duality
\[
h\of{\bar{y}^*, \bar{z}^*} \supseteq h\of{y^*, z^*} \supseteq
    p \oplus S\of{\bar{z}^*} = h\of{\bar{y}^*, \bar{z}^*},
\]
hence $h\of{\bar{y}^*, \bar{z}^*}$ is maximal in $H$.

(iii) Assume $h\of{\bar{y}^*, \bar{z}^*}$ is maximal in $H$. Then, there is
$\bar{\bar{y}}^* \in Y^*$ such that $Z \neq p \oplus S\of{\bar{z}^*} =
h\of{\bar{\bar{y}}^*, \bar{z}^*}$ according to the strong duality theorem. Weak duality
implies
\[
h\of{\bar{y}^*, \bar{z}^*} \supseteq p \oplus S\of{\bar{z}^*} = h\of{\bar{\bar{y}}^*,
\bar{z}^*},
\]
hence $h\of{\bar{y}^*, \bar{z}^*} = h\of{\bar{\bar{y}}^*, \bar{z}^*}$ by maximality of
$h\of{\bar{y}^*, \bar{z}^*}$. This proves $\of{\bar{y}^*, \bar{z}^*} \in \Delta$. \pend

\subsection{Strong duality via compactness}
\label{SubSecCompact}

In this subsection, we use a "dual" condition to establish strong duality.

\begin{theorem}
\label{ThmCompactDual}
Let $f \colon X \to \mathcal G\of{C}$ and $g \colon X \to \mathcal G\of{D}$ be closed. Assume that there are $z^*_0 \in C^+$ and $y^*_0 \in D^+$ such that the function
\[
x \mapsto \vp_{f, z^*_0}\of{x} + \vp_{g, y^*_0}\of{x} 
\]
has compact sublevel sets. Then, the value function $v\of{y} = \inf\cb{f\of{x} \mid y \in g\of{x}}$ is closed. If, additionally, $f$ and $g$  are convex, then strong duality holds, i.e. $v\of{0} = p = d$.
\end{theorem}

{\sc Proof.} Take a net $\of{y_\alpha, z_\alpha} \in \gr v$ which converges to $\of{y, z} \in Y \times Z$. Then, there is a net $x_\alpha$ and for each $\alpha$ there is a net $z_\beta \to 0$ such that
\[
z_\alpha + z_\beta \in f\of{x_\alpha}, \quad y_\alpha \in g\of{x_\alpha}.
\]
Since $f$ and $g$ have closed convex values, this implies (see \eqref{EqReScalarize})
\[
z^*_0\of{z_\alpha + z_\beta} \geq \vp_{f, z^*_0}\of{x_\alpha}, \quad y^*_0\of{y_\alpha} \geq \vp_{g, y^*_0}\of{x_\alpha}.
\]
Since $y_\alpha$, $z_\alpha$ and $z_\beta$ are convergent, there is $\alpha_0$ such that for all $\alpha$ exceeding $\alpha_0$ there is $\beta_0\of{\alpha}$ such that
\[
\vp_{f, z^*_0}\of{x_\alpha} + \vp_{g, y^*_0}\of{x_\alpha} \leq z^*_0\of{z_\alpha + z_\beta} + y^*_0\of{y_\alpha} \leq z^*_0\of{z} + y^*_0\of{y} + 1 
\]
for all $\beta$ exceeding $\beta_0\of{\alpha}$. Because of the compact sublevel set assumption, $x_\alpha$ has a subnet converging to some $x \in X$. Since $f$ and $g$ are closed, we have $z \in f\of{x}$ and $y \in g\of{x}$, hence $z \in v\of{y}$. This proves that $v$ is indeed closed. If $v\of{0} = Z$, then $d = p$ by weak duality. If the additional convexity assumption is satisfied and $v\of{0} \neq Z$, then $v$ is closed convex and $\of{\cl v}\of{0} = v\of{0} = v^{**}\of{0}$. Proposition \ref{PropDualOptimalValue} (b) implies strong duality. \pend

\medskip Note that the same proof works if one replaces the function $\vp_{f, z^*_0} + \vp_{g, y^*_0}$ in the compact sublevel set assumption by $\lambda_0\vp_{f, z^*_0} + \lambda_1\vp_{g, y^*_0}$ for real numbers $\lambda_0, \lambda_1 \geq 0$. The remarkable fact is that such an assumption for just one scalarization is enough the ensure strong duality. This is in the spirit of \cite[Theorem 4.4]{HeydeSchrage11R} and \cite[Theorem 3.9]{HamelSchrage12R}. On the other hand, primal attainment of the solution is a much harder to achieve property which will be discussed elsewhere.

%Newsubsection
\subsection{An example}
\label{SubSecExample}

We shall consider the following specialization: Let $p, q, N, M$ be positive integers and $X = \R^{pN}$, $Y = \R^M$, $Z = \R^q$ and $C = \R^q_+$, $D = \R^m_+$. Moreover, let the matrices $A_n \in \R^{M \times p}$, $n \in \cb{1,2, \ldots, N}$, and $b \in \R^M$ be given. The problem we are interested in is
\[
\tag{PSep} \mbox{minimize} \quad \sum_{n=1}^Nf_n\of{x^n} \quad \mbox{subject to} \quad  0 \in \sum_{m=1}^MA_nx^n - b + \R^M_+
\]
for functions $f_n \colon \R^p \to \G\of{\R^q, \R^q_+}$, $n = 1, \ldots, N$. The motivation for this problem with separated variables comes from utility maximization for vector-valued utility functions. One may compare \cite[Chapter 3]{DelbaenSchachermayer06} for a scalar version and also \cite[Section 4.3, Example 4(d)]{BorweinLewis06Book}.

The Lagrangian for this problem is
\begin{align*}
l\of{x, v, w} & =  \sum_{n=1}^Nf_n\of{x^n} + S_{\of{-v, w}}\of{\sum_{n=1}^N A_n x^n - b}\\
	& = \sum_{n=1}^Nf_n\of{x^n} + \sum_{n=1}^N S_{\of{-A^T_n v, w}}\of{x^n} + S_{\of{-v, w}}\of{-b}.
\end{align*}
Here, by a slight abuse of notation, the dual variables are $v \in \R^M$, $w \in \R^q_+\bs\cb{0}$. The dual objective becomes
\begin{align*}
h\of{v, w} &  =  \inf_{x \in \R^{pN}} l\of{x, v, w} \\
	& = S_{\of{v, w}}\of{b} + \inf_{x^1, \ldots, x^N \in \R^p} \sum_{n=1}^N\sqb{f_n\of{x^n} + S_{\of{-A^T_n v, w}}\of{x^n}} \\
	& = S_{\of{v, w}}\of{b} + \sum_{n=1}^N\inf_{x^n \in \R^p} \sqb{f_n\of{x^n} + S_{\of{A^T_n v, w}}\of{-x^n}} \\
	& = S_{\of{v, w}}\of{b} + \sum_{n=1}^N -f^*_n\of{A^T_n v, w}.
\end{align*}
If one defines a function $U \colon \R^M \to \mathcal G\of{\R^q, \R^q_+}$ by
\[
U\of{y} = \inf_{x \in \R^{pN}}\cb{ \sum_{n=1}^Nf_n\of{x^n} \,\vert\, 0 \in \sum_{m=1}^MA_n x^n - y + \R^M_+}
\]
and another one $-V \colon \R^M \times \R^q_+\bs\cb{0}$ by
\[
-V\of{v, w} =  \sum_{n=1}^N -f^*_n\of{A^T_n v, w},
\]
then, in accordance with well-known formulas from scalar utility optimization (see \cite[Chapter 3]{DelbaenSchachermayer06}) strong duality is nothing else than
\[
U\of{b} = \sqb{-V\of{\cdot, \cdot}}^*\of{-b}.
\]

%Newsubsection

\subsection{Set-valued duality in vector optimization}

In \cite[Theorem 3.31, 3.32]{Loehne11Book}, Lagrange duality theorems for infimal-set-valued functions are given under the assumption that the interior of the ordering cone $C$ is non-empty. Here, we show that comparable results can be obtained under a much weaker assumption. Along the way, the relationship between the Fenchel conjugates introduced by  \cite{Loehne05Diss,Loehne05} and those by \cite{Hamel09} will be clarified.

Let us assume that there is an element $z_0 \in C\bs\cb{0}$ such that
\[
\forall z^* \in C^+\bs\cb{0} \colon z^*\of{z} > 0.
\]
In this case, the set $B\of{z_0} = \cb{z^* \in C^+ \colon z^*\of{z_0} = 1}$ is a base of
$C^+$ with $0 \not\in \cl B\of{z_0}$. That is, for each $z^* \in C^+\bs\cb{0}$ there is a
unique representation $z^* = tz^*_0$ with $t>0$ and $z^*_0 \in B\of{z_0}$. Compare
\cite{GoeRiaTamZal03}, Definition 2.1.14, Theorem 2.1.15 and 2.2.12 applied to $C^+$
instead of $C$. Clearly, a pointed closed convex cone with non-empty interior has a base, and conversely, the cone $L^2_+$ has an empty interior, but a base is generated by the constant 1 function.

The very definition of the functions $S_{\of{x^*, z^*}}$ gives
\[
\cb{S_{\of{x^*, z^*}} \mid x^* \in X^*, \; z^* \in C^+\bs\cb{0}} =
    \cb{S_{\of{x^*, z^*}} \mid x^* \in X^*, \; z^* \in B\of{z_0}}.
\]
Therefore, it is sufficient to run an intersection like in the definition of $d$ in
Section \ref{SubsecDual} over $y^* \in Y^*$ and $z^* \in B\of{z_0}$. Moreover, one easily checks (see also Proposition 6 (iv) in \cite{Hamel09}) for $z^* \in B\of{z_0}$
\[
\forall x \in X \colon S_{\of{x^*, z^*}}\of{x} = \cb{x^*\of{x}z_0} + S\of{z^*}.
\]
Thus, the conjugate of a function $f \colon X \to \P\of{C}$ can be written as
\[
-f^*\of{x^*, z^*} =
    \cl\bigcup_{x \in X} \sqb{f\of{x} - x^*\of{x}z_0 + S\of{z^*}} =
    \cl\bigcup_{x \in X} \sqb{f\of{x} - x^*\of{x}z_0} \oplus S\of{z^*}.
\]
The part which does not depend on $z^*$ (remember $z_0$ defines a base of $C^+$ and is
the same for all $z^* \in C^+\bs\cb{0}$) has been used in \cite{Loehne05},
\cite{LoehneTammer07} for a definition of another set-valued conjugate, namely
\[
-f^*_{z_0}\of{x^*} = \cl \bigcup_{x \in X} \sqb{f\of{x} - x^*\of{x}z_0}.
\]
If the cone $C$ has non-empty interior and $z_0 \in \Int C$, then the space $\mathcal{F}$
of upper closed sets as defined in \cite{Loehne11Book} coincides with $\mathcal F\of{C}$. In this case, one can derive an infimal set version of the Lagrange duality theorem (Theorem \ref{ThmStrongDuality}) in the same way as
Theorem 3.32 is derived from Theorem 3.26 in \cite{Loehne11Book}. Note that under
convexity assumptions the functions map indeed into $\mathcal G\of{C}$ which is a subset of
$\mathcal F\of{C}$. We omit the details since they can be found in \cite{Loehne11Book}. Finally, if $Z = \R$, $C = \R_+$, then $C^+ = \R_+$, and $\cb{1}$ is a base of $C^+$. With this simple device one obtains a scalar version of Theorem \ref{ThmStrongDuality}, for example, Theorem 4.3.7. in \cite{BorweinLewis06Book}.

%New section
\section{Canonical extensions and saddle points of Lagrangians}
\label{SecSaddle}

%New subsection
\subsection{Saddle points for functions with values in complete lattices}

With a bi-variable function $l \colon X \times W \to \mathcal L$, one can associate two
optimization problems, namely
\begin{align*} \tag{Pl}
 & \mbox{minimize} \quad p\of{x} \quad \mbox{over} \quad x \in X \quad \mbox{and} \\
 \tag{Dl} & \mbox{maximize} \quad d\of{w} \quad \mbox{over} \quad w \in W
\end{align*}
where
\begin{align*}
p\of{x} & = \hat{L}\of{\cb{x}, W} = \sup_{w \in W}l\of{x, w}, \\
d\of{w} & = \check{L}\of{X, \cb{w}} = \inf_{x \in X}l\of{x, w}.
\end{align*}

The weak duality relation
\[
\forall x \in X, \; \forall w \in W \colon d\of{w} \leq p\of{x}
\]
follows immediately from \eqref{EqLowerLessUpper}. Indeed, for all $x \in X$ and all $w
\in W$,
\[
d\of{w} = \hat{L}\of{X, \cb{w}} \leq
    \hat{L}\of{X, W}
    \stackrel{\eqref{EqLowerLessUpper}}{\leq} \check{L}\of{X, W} \leq
    \check{L}\of{\cb{x}, W} = p\of{x}.
\]

In the next definition, saddle points for the function $l$ are introduced.

\begin{definition}
\label{DefSaddlePoints} Let $l \colon X \times V \to \mathcal L$ be a function. A pair
$\of{\bar X, \bar V} \in \mathcal P\of{X} \times \mathcal P\of{V}$ is called a {\em
saddle point} of $l$ if the following conditions are satisfied:
\\
(a) $\emptyset \neq p\sqb{\bar X} \subseteq \Min p\sqb{X}$ and $\emptyset \neq d\sqb{\bar
V} \subseteq \Max d\sqb{V}$.
\\
(b) For all $U \in \mathcal P\of{X}$ and all $W \in \mathcal P\of{V}$,
\[
\check{L}\of{\bar X, W} \leq \check{L}\of{\bar X, \bar V}
 = \hat{L}\of{\bar X, \bar V} \leq \hat{L}\of{U,
\bar V}.
\]
A saddle point is called {\em full} if the two inclusions in (a) are equations.
\end{definition}

Condition (b) generalizes the saddle point condition known from scalar optimization. Note
that (b) also includes the statement that the lower and the upper canonical extension of
$l$ coincide at the saddle point. In the scalar case, (a) and (b) are equivalent -- which
is no longer true for functions with values in complete lattices.

\begin{lemma}
\label{LemInfSupSaddle} For a pair $\of{\bar X, \bar V} \in \mathcal P\of{X} \times
\mathcal P\of{V}$, statement (b) in Definition \ref{DefSaddlePoints} is equivalent to
\begin{equation}
\label{EqInfSupSaddle}
 \sup_{v \in \bar V} d\of{v} = \inf_{x \in \bar X} p\of{x}.
\end{equation}
\end{lemma}

{\sc Proof.} Assume (b) of Definition \ref{DefSaddlePoints}. With $W = V$ and $U = X$ we
obtain
\[
\inf_{x \in \bar X} p\of{x} = \check{L}\of{\bar X, V}
    \leq \hat{L}\of{X, \bar V} = \sup_{v \in \bar V} d\of{v}
\]
which produces "$\leq$" in \eqref{EqInfSupSaddle} according to the definitions of
$\check{L}$, $\hat{L}$, $p$ and $d$. On the other hand, from the
weak duality relation \eqref{EqLowerLessUpper} we obtain
\[
\sup_{v \in \bar V} d\of{v} =
    \hat{L}\of{X, \bar V} \leq \hat{L}\of{\bar X, \bar V}
    \leq \check{L}\of{\bar X, \bar V} \leq \check{L}\of{\bar X, V}
    = \inf_{x \in \bar X} p\of{x}.
\]

Conversely, assume \eqref{EqInfSupSaddle}. Then, for $U \in \mathcal P\of{X}$, $W \in
\mathcal P\of{V}$,
\[
\inf_{x \in \bar X} \sup_{v \in W} l\of{x, v} \leq
 \inf_{x \in \bar X} \sup_{v \in V} l\of{x, v}
= \sup_{v \in \bar V} \inf_{x \in X} l\of{x, v} \leq
 \sup_{v \in \bar V} \inf_{x \in U} l\of{x, v}
\]
which produces
\[
\check{L}\of{\bar X, W} \leq \hat{L}\of{U, \bar V}.
\]
From this, we obtain
\begin{align*}
\forall U \in \mathcal P\of{X} & \colon \check{L}\of{\bar U, \bar V} \leq
 \hat{L}\of{U, \bar V} \\
\forall W \in \mathcal P\of{V} & \colon
 \check{L}\of{\bar X, W} \leq \hat{L}\of{\bar X, \bar V}
\end{align*}
as well as
\[
\check{L}\of{\bar X, \bar V} \leq
 \hat{L}\of{\bar X, \bar V}.
\]
The weak duality relation \eqref{EqLowerLessUpper} gives
\[
\hat{L}\of{\bar X, \bar V} \leq
 \check{L}\of{\bar X, \bar V},
\]
and the last four relations together imply (b) of Definition \ref{DefSaddlePoints}. \pend

\begin{theorem}
\label{ThmPDSolSaddle} The following statements are equivalent:
\\
(a) $\bar X \in \mathcal P\of{X}$ is a (full) solution of (Pl), $\bar V \in \mathcal
P\of{W}$ is a (full) solution of (Dl), and strong duality holds.
\\
(b) $\of{\bar X, \bar V}$ is a (full) saddle point of $l \colon X \times V \to \mathcal
L$.
\end{theorem}

{\sc Proof.} (a) $\Rightarrow$ (b): Since $\bar X$ is a solution of of (Pl) and $\bar V$
of (Dl), (a) of Definition \ref{DefSaddlePoints} is immediate. Strong duality and the
attainment of the infimum/supremum yield
\[
 \sup_{v \in \bar V} d\of{v} =  \sup_{v \in V} d\of{v}
    =  \inf_{x \in X} p\of{x} = \inf_{x \in \bar X} p\of{x},
\]
hence Lemma \ref{LemInfSupSaddle} produces (b) of Definition \ref{DefSaddlePoints}.

(b) $\Rightarrow$ (a): From (a) of Definition \ref{DefSaddlePoints} we obtain that $\bar
X$ is a minimizer of (Pl) and $\bar V$ a maximizer of (Dl). Lemma \ref{LemInfSupSaddle}
and weak duality yield
\[
\inf_{x \in \bar X} p\of{x} = \sup_{v \in \bar V} d\of{v}
    \leq \sup_{v \in V} d\of{v} \leq \inf_{x \in X} p\of{x} \leq \inf_{x \in \bar X}
    p\of{x},
\]
hence $\bar X$ is an infimizer of (Pl), $\bar V$ a supremizer of (Dl) and strong duality
holds. \pend

%New subsection
\subsection{Saddle points of the Lagrangian and strong duality}

In this subsection, we specify the setting of the previous one to the set-valued
optimization problems introduced in Section \ref{SecProblem}. In particular, with
$\of{\mathcal L, \leq} = \of{\G, \supseteq}$ and $W = Y^* \times C^+\bs\cb{0}$ we
can define the Lagrangian of problem (P) as the function
\[
\of{x, y^*, z^*} \mapsto l\of{x, y^*, z^*} =  f\of{x} \oplus \bigcup_{y \in
g\of{x}}S_{\of{y^*, z^*}}\of{y}
\]
which depends on the primal variable $x \in X$ and the pair of dual variables $\of{y^*,
z^*} \in Y^* \times C^+\bs\cb{0}$.

\begin{remark}
Of course, if $l$ is the Lagrangian of problem (P), then, under the assumption of
Proposition \ref{PropReconPrimal}, (Pl) coincides with (P), and (Dl) coincides with (D)
since, by definition, $d = h$ in this case.
\end{remark}

The following theorem gives the link between saddle points of the Lagrangian and strong
duality.

\begin{theorem}
\label{ThmStrongDualitySaddle} Let the assumptions of Theorem \ref{ThmStrongDuality} be
satisfied, and let $\bar U \in \mathcal P\of{X}$ be a solution of (P). Then, there exists
$\bar \Delta \in \mathcal P\of{Y^* \times C^+\bs\cb{0}}$ such that $\of{\bar U, \bar
\Delta}$ is a saddle point of the Lagrangian $l$.
\end{theorem}

{\sc Proof.} Follows from Theorem \ref{ThmPDSolSaddle} and Theorem
\ref{ThmStrongDuality}. \pend

%%%Newsection
\section{Appendix}

The following definition is taken from \cite{Hamel05Habil} where references and more
material about structural properties of conlinear spaces can be found.

\begin{definition}
\label{DefConlinearSpace} A nonempty set $W$ together with two algebraic operations $+
\colon W \times W \to W$ and $\cdot \colon \R_+ \times W \to W$ is called a conlinear
space provided that
\\
(C1) $\of{W, +}$ is a commutative monoid with neutral element $\theta$,
\\
(C2) (i) $\forall w_1, w_2 \in W$, $\forall r \in \R_+$: $r \cdot \of{w_1 + w_2} = r
\cdot w_1 + r \cdot w_2$, (ii) $\forall w \in W$, $\forall r, s \in \R_+$: $s \cdot \of{r
\cdot w} = \of{rs} \cdot w$, (iii) $\forall w \in W$: $1 \cdot w = w$, (iv) $0 \cdot
\theta = \theta$.

An element $w \in W$ is called a convex element of the conlinear space $W$ if
\[
\forall s, t \geq 0 \colon \of{s+t} \cdot w = s \cdot w + t \cdot w.
\]

A conlinear space $\of{W, +, \cdot}$ together with a partial order $\preceq$ on $W$ (a
reflexive, antisymmetric, transitive relation) is called ordered conlinear space provided
that (iv) $w, w_1, w_2 \in W$, $w_1 \preceq w_2$ imply $w_1 + w \preceq w_2 + w$, (v)
$w_1, w_2 \in W$, $w_1 \preceq w_2$, $r \in \R_+$ imply $r \cdot w_1 \preceq r\cdot w_2$.

A non-empty subset $V \subseteq W$ of the conlinear space $\of{W, +, \cdot}$ is called a conlinear subspace of $W$ if (vi) $v_1, v_2 \in V$ implies $v_1 + v_2 \in V$ and (vii) $v \in V$ and $t \geq 0$ imply $t \cdot v \in V$.
\end{definition}

It can easily be checked that a conlinear subspace of a conlinear space again is a conlinear space.

\bibliographystyle{plain}
\bibliography{HamelLoehneBIB} 

\end{document}